\documentclass[11pt]{article}
\usepackage{latexsym,amsfonts,amsmath,amsthm,amssymb, epsfig, bm}
\usepackage[english]{babel}
\usepackage[margin=1.5in]{geometry}
\usepackage{color}

\newcommand{\vek}[1]{\boldsymbol{#1}}
\newcommand{\norm}[2]{\left\|\left.{#1}\right|{#2}\right\|}

\newcommand{\p}{p(\cdot)}
\newcommand{\q}{q(\cdot)}
\newcommand{\s}{{s(\cdot)}}

\newcommand{\SSn}{\mathcal{S}'(\rn)}

\newcommand{\ellqp}{{\ell_{q(\cdot)}(L_{p(\cdot)}(\rn))}}
\newcommand{\ellpq}{{L_{p(\cdot)}(\ell_{q(\cdot)}(\rn))}}
\newcommand{\Bwpqpunkt}{B^{\bm{w}}_{\p,\q}(\rn)}
\newcommand{\Fwpqpunkt}{F^{\bm{w}}_{\p,\q}(\rn)}
\newcommand{\Bspq}{B^{s}_{p,q}(\rn)}
\newcommand{\Fspq}{F^{s}_{p,q}(\rn)}
\newcommand{\Bwpq}{B^{\bm{w}}_{p,q}(\rn)}
\newcommand{\Fwpq}{F^{\bm{w}}_{p,q}(\rn)}
\newcommand{\Bspqpunkt}{B^{s(\cdot)}_{\p,\q}(\rn)}
\newcommand{\Fspqpunkt}{F^{s(\cdot)}_{\p,\q}(\rn)}
\newcommand{\bspqpunkt}{b^{s(\cdot)}_{\p,\q}(\rn)}
\newcommand{\fspqpunkt}{f^{s(\cdot)}_{\p,\q}(\rn)}

\newcommand{\Plog}{\mathcal{P}^{\log}(\rn)}

\newcommand{\supp}{\mathop{\mathrm{supp}\,}\nolimits}

\newcommand{\real}{\mathbb{R}}
\newcommand{\rn}{{\mathbb{R}^n}}
\newcommand{\nat}{\mathbb{N}}

\newcommand{\zn}{{\mathbb{Z}^n}}

\newtheorem{definition}{Definition }[section]
\newtheorem{theorem}[definition]{Theorem}

\newtheorem{corollary}[definition]{Corollary}
\newtheorem{lemma}[definition]{Lemma}

\newtheorem{remark}[definition]{Remark}

\newcommand{\beq}{\begin{equation}}
\newcommand{\eeq}{\end{equation}}

\newcommand{\bpr}{\begin{proof}}
\newcommand{\epr}{\nopagebreak \hfill $\square$\\ \end{proof}}

\title{Franke-Jawerth embeddings for Besov and Triebel-Lizorkin spaces with variable exponents}
\author{Helena F. Gon\c{c}alves\footnote{{\tt helena.goncalves@mathematik.tu-chemnitz.de}; the author was supported by the German science foundation (DFG) within the project KE 1847/1-1. }, Henning Kempka\footnote{{\tt henning.kempka@mathematik.tu-chemnitz.de}; the author was supported by the German science foundation (DFG) within the project KE 1847/1-1.},
Jan Vyb\'iral\footnote{Department of Mathematical Analysis, Charles University, Sokolovsk\'a 83, 186 00, Prague 8, Czech Republic, {\tt vybiral@karlin.mff.cuni.cz};
this author was supported by the ERC CZ grant LL1203 of the Czech Ministry of Education and by the Neuron Fund for Support of Science.}}

\begin{document}
\maketitle
\begin{abstract}
The classical Jawerth and Franke embeddings 
\beq
F^{s_0}_{p_0,q}(\rn)\hookrightarrow B^{s_1}_{p_1,p_0}(\rn) \quad \mbox{and} \quad B^{s_0}_{p_0,p_1}(\rn)\hookrightarrow F^{s_1}_{p_1,q}(\rn) \nonumber
\eeq
are versions of Sobolev embedding between the scales of Besov and Triebel-Lizorkin function spaces for $s_0>s_1$ and $\displaystyle s_0-\frac{n}{p_0} = s_1 -\frac{n}{p_1}.$
We prove Jawerth and Franke embeddings for the scales of Besov and Triebel-Lizorkin spaces with all exponents variable
\beq
F^{s_0(\cdot)}_{p_0(\cdot),\q}\hookrightarrow B^{s_1(\cdot)}_{p_1(\cdot),p_0(\cdot)} \quad \mbox{and} \quad B^{s_0(\cdot)}_{p_0(\cdot),p_1(\cdot)}\hookrightarrow F^{s_1(\cdot)}_{p_1(\cdot),\q}, \nonumber
\eeq
respectively, if $\inf_{x\in\mathbb{R}^n}(s_0(x)-s_1(x))>0$ and
\beq
s_0(x) -\frac{n}{p_0(x)} = s_1(x) -\frac{n}{p_1(x)}, \quad x \in \rn. \nonumber
\eeq
We work exclusively
with the associated sequence spaces $\bspqpunkt$ and $\fspqpunkt$, which is justified by well known decomposition techniques. 
We give also a different proof of the Franke embedding in the constant exponent case which avoids duality arguments and interpolation.

Our results hold also for 2-microlocal function spaces $\Bwpqpunkt$ and $\Fwpqpunkt$ which unify the smoothness scales of spaces of variable smoothness and generalized smoothness spaces. 

\end{abstract}



\section{Introduction}

\indent

Spaces of variable integrability, also known as variable exponent function spaces $L_{\p}(\rn)$, can be traced back to Orlicz \cite{Or31} 1931, but the modern development started with the papers \cite{KR91} of Kov\'a\v{c}ik and R\'akosn\'ik as well as
\cite{ER00} of Edmunds and R\'akosn\'\i k and \cite{Die04} of Diening. The spaces $L_{\p}(\rn)$ have interesting applications in fluid dynamics, namely in the theory of electrorheological fluids \cite{Ruzicka}, where $p(\cdot)$ is a function of the electric field. Further, these variable function spaces were used in image processing, PDEs and variational calculus, see the introduction of \cite{DHR09}. For an overview we refer to \cite{DHHR}.

Sobolev and Besov spaces with variable smoothness but fixed integrability have been introduced in the late 60's and early 70's in the works of Unterberger \cite{Unterberger}, Vi\v{s}ik and Eskin \cite{VisikEskin}, Unterberger and Bokobza \cite{UnterbergerBok} and in the work of Beauzamy \cite{Beauzamy}. Leopold studied in \cite{Leopold91} Besov spaces where the smoothness is determined by a symbol $a(x,\xi)$ of a certain class of hypoelliptic pseudodifferential operators. In the special case $a(x,\xi)=(1+|\xi|^2)^{\sigma(x)/2}$ these spaces coincide with spaces of variable smoothness $B^{\sigma(x)}_{p,p}(\rn)$.

A more general approach to spaces of variable smoothness are the so-called 2-microlocal function spaces $\Bwpq$ and $\Fwpq$. The smoothness in these scales gets measured by a weight sequence $\bm{w}=(w_j)_{j\in \nat_0}$. Besov spaces with such weight sequences appeared first in the works of Peetre \cite{PeetreArt} and Bony \cite{Bony}. Establishing a wavelet characterization for 2-microlocal H\"older-Zygmund spaces in \cite{Jaf91} it turned out that 2-microlocal spaces are well adapted in connection to regularity properties of functions (\cite{jm96},\cite{Meyer97},\cite{LevySeu04}). Spaces of variable smoothness are a special case of 2-microlocal function spaces and in \cite{LevySeu03} and \cite{Bes1} characterizations by differences have been given for certain classes of them.\\

Even in the case of constant exponents the integrability exponents $p,q\in(0,\infty]$ and the smoothness parameter $s$ inherit a quite interesting interplay regarding embeddings and special cases, see \cite{Triebel2} and \cite{Triebel3}. 

If one considers at first Triebel-Lizorkin spaces $F^s_{\p,q}(\rn)$ where only the integrability parameter $\p$ is chosen to  be variable, then the other exponents must be chosen variable as well. This can already be seen by the Sobolev embedding from \cite{Vyb09}
\begin{align*}
	F^{s_0(\cdot)}_{p_0(\cdot),q}(\rn)\hookrightarrow F^{s_1(\cdot)}_{p_1(\cdot),q}(\rn)
\end{align*}
under the usual condition, but now pointwise,
\begin{align*}
	s_0(x)-\frac{n}{p_0(x)}=s_1(x)-\frac{n}{p_1(x)},\quad x\in\rn.
\end{align*} 
Now also the smoothness parameter $s(\cdot)$ should be chosen variable. That also the third index $\q$ should be variable can be seen by the following trace theorem. It was obtained by Diening, H{\"a}st{\"o} and Roudenko in \cite{DHR09} where variable smoothness and integrability were for the first time combined in one approach. They defined Triebel-Lizorkin spaces $F^{s(\cdot)}_{\p, \q}(\rn)$ and considered the trace theorem on $\mathbb{R}^{n-1}$. Here the usual result holds in a variable analogue
\[
\mathrm{Tr}\  F^{s(\cdot)}_{\p,\q}(\rn)={F}^{s(\cdot)-\frac{1}{\p}}_{\p,\p}(\mathbb{R}^{n-1}), \mbox{ with } s(\cdot) - \frac{1}{\p} > (n-1)\max\left(\frac{1}{\p}-1, 0\right),
\]
(Theorem~3.13 in \cite{DHR09}) and we see the necessity of taking $s$ and $q$ variable if $p$ is not constant.

For the Besov spaces it is non-trivial to have also the parameter $q$ as a variable one. Almeida and H{\"a}st{\"o} were able to introduce in \cite{AH10} Besov spaces $B_{\p,\q}^{s(\cdot)}(\rn)$ with all three indices variable and proved the Sobolev and other usual embeddings in this scale. These spaces need to be defined by using another modular which already uses the variable structure on $\q$, see Section \ref{sec:notation}.\\
Interestingly, these variable Besov spaces fit very well to the constant exponent case theory which can be seen by the embedding
\begin{align*}
B^{s(\cdot)}_{\p,\min(\p,\q)}(\rn)\hookrightarrow F^{s(\cdot)}_{\p,\q}(\rn)\hookrightarrow B^{s(\cdot)}_{\p,\max(\p,\q)}(\rn).
\end{align*}
On the other hand, in \cite{KemVybnorm} it has been shown that the triangle inequality in $B^{s(\cdot)}_{\p,\q}(\rn)$ is in general not true for exponents with $\min(\p,\q)\geq1$. This in sharp contrast to the case of Triebel-Lizorkin spaces $F^{\s}_{\p,\q}(\rn)$ and to the constant exponent spaces $B^s_{p,q}(\rn)$ and $F^s_{p,q}(\rn)$, which are always normed spaces if $\min(\p,\q)\geq1$, or $\min(p,q)\geq1$,  respectively. 

For the full variable spaces $B^{s(\cdot)}_{\p,\q}(\rn)$ and $F^{s(\cdot)}_{\p,\q}(\rn)$ different characterizations of the spaces as decompositions via atoms, molecules, local means and ball means of differences (see \cite{DHR09}, \cite{Dri12}, \cite{KemV12}) have been shown. Further, there also exist results on the extension operator from halfspaces in \cite{Noi}. 

Furthermore, also for the more general scale of 2-microlocal function spaces $B^{\vek{w}}_{\p,\q}(\rn)$ and $F^{\vek{w}}_{\p,\q}(\rn)$ (see Section \ref{sec:2ml} for details) all the above mentioned characterizations have been obtained and there exist results on traces \cite{GMN14}, pointwise multipliers \cite{GK15} and Fourier multipliers \cite{AlmCaet}.

Regarding Franke-Jawerth embeddings, they go back to Jawerth in \cite{Jaw77} and Franke in \cite{Fra86}. Using interpolation techniques {and duality}, the authors proved the following.
\begin{theorem}\label{thm:FJ}
 Let $-\infty < s_1<s_0<\infty$, $0<p_0<p_1\leq \infty$ and $0<q\leq \infty$ with
 $$s_0-\frac{n}{p_0} = s_1-\frac{n}{p_1}.$$
 \begin{list}{}{\labelwidth1.3em\leftmargin2em}
    \item[{\upshape (i)\hfill}] Then \begin{equation}\label{eq:FJ1}F^{s_0}_{p_0,q}(\rn) \hookrightarrow B^{s_1}_{p_1,p_0}(\rn).\end{equation}
    \item[{\upshape (ii)\hfill}] If $p_1<\infty$, then \begin{equation}\label{eq:FJ2} B^{s_0}_{p_0,p_1}(\rn) \hookrightarrow F^{s_1}_{p_1,q}(\rn).\end{equation}
 \end{list}
\end{theorem}
The surprising effect in \eqref{eq:FJ1} and \eqref{eq:FJ2} is that (unlike in the case of Sobolev embeddings)
no conditions on $q$ are necessary. This shows that Jawerth and Franke embeddings
exploit the fine properties of Besov and Triebel-Lizorkin spaces and exhibit an interesting interplay between these two scales of function spaces.

Later Vyb{\'{\i}}ral in \cite{Vyb08} gave a new proof of Theorem \ref{thm:FJ}. The author transferred the problem to the corresponding sequence spaces and, instead of interpolation, he used the technique of non-increasing rearrangements as well as duality.
The developed technique via sequence spaces was also used in \cite{HaroskeLeszek} to obtain the Franke-Jawerth embeddings in the Morrey space versions of Besov and Triebel-Lizorkin spaces and in \cite{HansenVybiral} to get these embeddings for spaces with dominating mixed smoothness.

Our aim is to extend these results to the scale of Besov and Triebel-Lizorkin spaces with variable smoothness and integrability $\Bspqpunkt$ and $\Fspqpunkt$, obtaining in this way a fine connection
between those scales of function spaces. One can observe then that the somehow artificial definition of $\Bspqpunkt$ seems to be well chosen.\\

The paper is organized as follows. We introduce in Section \ref{sec:notation} the necessary notation and definitions which are needed afterwards. Furthermore, we also state same known theorems for the spaces with variable exponents. In Section \ref{sec:FrankeKonstant} we present another proof for the Franke embedding in the constant exponent case.
The novelty of the technique here is that we totally avoid the use of interpolation and duality arguments. With the help of this result in the constant exponent case we state and prove the Jawerth and Franke embeddings in the scales of Besov $\Bspqpunkt$ and Triebel-Lizorkin spaces $\Fspqpunkt$ with variable exponents in Section \ref{sec:FrankeJawerth}.
In Section \ref{sec:2ml} we transfer our results to 2-microlocal function spaces with variable exponents.
Finally, in the last section we pose some open problems. 

\section{Notation and definitions}\label{sec:notation}
\indent

We shall adopt the following general notation: $\nat$ denotes the set of all natural numbers, $\nat_0=\mathbb N\cup\{0\}$, $\mathbb{Z}$ denotes the set of integers, $\rn$ for $n\in\nat$ denotes the $n$-dimensional real Euclidean space with $|x|$, for $x\in\rn$, denoting the Euclidean norm of $x$.

For $q\in (0,\infty]$, $\ell_q$ stands for the linear space of all complex sequences ${a}=(a_j)_{j\in \nat_0}$ endowed with the quasi-norm
$$
\Vert {a} \mid \ell_q\Vert = \Big(\sum_{j=0}^{\infty} |a_j|^q \Big)^{1/q},
$$
with the usual modification if $q=\infty$. By $c$, $C$, etc. we denote positive constants independent of appropriate quantities. For two non-negative expressions ({\it i.e.}, functions or functionals) ${\mathcal  A}$, ${\mathcal  B}$, the symbol ${\mathcal A}\lesssim {\mathcal  B}$ (or ${\mathcal A}\gtrsim {\mathcal  B}$) means that $ {\mathcal A}\leq c\, {\mathcal  B}$ (or $c\,{\mathcal A}\geq {\mathcal B}$), for some $c>0$. If ${\mathcal  A}\lesssim {\mathcal  B}$ and ${\mathcal A}\gtrsim{\mathcal  B}$, we write ${\mathcal  A}\sim {\mathcal B}$ and say that ${\mathcal  A}$ and ${\mathcal  B}$ are equivalent.\\

Before introducing the function spaces under consideration we still need to recall some notation. By $\mathcal{S}(\rn)$ we denote the Schwartz space of all complex-valued rapidly decreasing infinitely differentiable functions on $\rn$ and by $\mathcal{S}'(\rn)$ its dual space of all tempered distributions on $\rn$. For $f\in \mathcal{S}'(\rn)$ we denote by $\widehat{f}$ the Fourier transform of $f$ and by $f^{\vee}$ the inverse Fourier transform of $f$. \\

Let $\varphi_0\in\mathcal{S}(\rn)$ be such that
\begin{equation}  \label{phi}
  \varphi_0(x)=1 \quad \mbox{if}\quad |x|\leq 1 \quad \mbox{and} \quad \supp \varphi_0 \subset \{x\in\rn: |x|\leq 2\}.
\end{equation}
Now define $\varphi(x):=\varphi_0(x)-\varphi_0(2x)$ and set $\varphi_j(x):=\varphi(2^{-j}x)$ for all $j\in\nat$. Then the sequence $(\varphi_j)_{j\in\nat_0}$ forms a smooth dyadic partition of unity.\\

By $\mathcal{P}(\rn)$ we denote the class of exponents, which are measurable functions $p:\rn\rightarrow (c,\infty]$ for some $c>0$. Let $p \in \mathcal{P}(\rn)$. Then, $p^+:=\text{ess-sup}_{x\in \rn}p(x)$, $p^-:=\text{ess-inf}_{x\in \rn}p(x)$ and $L_{p(\cdot)}(\rn)$ is the variable exponent Lebesgue space, which consists of all measurable functions $f$ such that for some $\lambda>0$ the modular $\varrho_{L_{p(\cdot)}(\rn)}(f/\lambda)$ is finite, where
$$
\displaystyle \varrho_{L_{p(\cdot)}(\rn)}(f):=\int_{\mathbb{R}^n_{0 }} |f(x)|^{p(x)}\, dx + \text{ess-sup}_{x\in \mathbb{R}^n_{\infty}} |f(x)|.
$$
Here $\mathbb{R}^n_{\infty}$ denotes the subset of $\rn$ where $p(x)=\infty$ and $\mathbb{R}^n_{0}= \rn \setminus \mathbb{R}^n_{\infty}$. The Luxemburg norm of a function $f\in L_{p(\cdot)}(\rn)$ is given by
$$
\|f \mid L_{p(\cdot)}(\rn)\|:=\inf\left\{\lambda>0:\varrho_{L_{p(\cdot)}(\rn)}\left(\frac{f}{\lambda}\right)\leq 1 \right\}.
$$

In order to define the mixed spaces $\ell_{\q}(L_{\p})$, we need to define another modular. For $p,q \in \mathcal{P}(\rn)$ and a sequence $(f_{\nu})_{\nu \in \nat_0}$ of complex-valued Lebesgue measurable functions on $\rn$, we define
\beq \label{modular_mixed}
  \varrho_{\ell_{\q}(L_{\p})}(f_{\nu}) = \sum_{\nu=0}^{\infty} \inf\left\{\lambda_{\nu}>0 : \varrho_{\p}\left(\frac{f_{\nu}}{\lambda_{\nu}^{1/\q}} \right)\leq 1 \right\}.
\eeq
If $q^+ <\infty$, then we can replace \eqref{modular_mixed} by the simpler expression
\beq \label{modular_mixed_norm}
  \varrho_{\ell_{\q}(L_{\p})}(f_{\nu}) = \sum_{\nu=0}^{\infty} \Big\| |f_{\nu}|^{\q} \mid L_{\frac{\p}{\q}} (\rn)\Big\|.
\eeq
The (quasi-)norm in the $\ell_{\q}(L_{\p})$ spaces is defined as usual by
\beq
  \| f_{\nu} \mid \ell_{\q}(L_{\p}(\rn)) \| = \inf \left\{ \mu>0 : \varrho_{\ell_{\q}(L_{\p})}\left(\frac{f_{\nu}}{\mu}\right) \leq 1\right\}.
\eeq

For the sake of completeness, we state also the definition of the space $L_{\p}(\ell_{\q})$. At first, one just takes the norm $\ell_{\q}$ of $(f_{\nu}(x))_{\nu \in \nat_0}$ for every $x \in \rn$ and then the $L_{\p}$-norm with respect to $x \in \rn$, i.e.
\beq
  \| f_{\nu} \mid L_{\p}(\ell_{\q}(\rn)) \| = \left\| \left(\sum_{\nu=0}^{\infty} |f_{\nu}(x)|^{q(x)} \right)^{1/q(x)} \mid L_{\p}(\rn)\right\|. \nonumber
\eeq \\

The following regularity classes for the exponents are necessary to make the definition of the spaces independent on the chosen decomposition of unity.

\begin{definition}\label{log-holder}
  Let $g\in C(\rn)$. We say that $g$ is locally log-H\"older continuous, abbreviated $g\in C_{loc}^{\log}(\rn)$, if there exists $c_{\log}(g)>0$ such that
  \beq \label{loc-log-holder}
      |g(x)-g(y)|\leq \frac{c_{\log}(g)}{\log(e+1/|x-y|)} \quad \text{for all}\;\; x,y\in\rn.
   \eeq
  We say that $g$ is globally log-H\"older continuous, abbreviated $g\in C^{\log}(\rn)$, if $g$ is locally log-H\"older continuous and there exists $g_{\infty}\in\real$ such that
    \beq
      |g(x)-g_{\infty}|\leq \frac{c_{\log}}{\log(e+|x|)} \quad \text{for all}\;\; x\in\rn.
    \eeq
\end{definition}
\vspace{0.5cm}
We use the notation $p\in \mathcal{P}^{\log}(\rn)$ if $p\in \mathcal{P}(\rn)$ and $1/p\in C^{\log}(\rn)$. It was proved in \cite{DHHMS09} that the maximal operator $\mathcal M$ is bounded in $L_{p(\cdot)}(\rn)$ provided that $p\in \mathcal{P}^{\log}(\rn)$ and $1<p^-\leq p^+\leq \infty$.\\

We recall the definition of the spaces $\Bspqpunkt$ and $\Fspqpunkt$, as given in \cite{DHR09} and \cite{AH10}.

\begin{definition}\label{def:BFpunkt}
Let $p,q\in \Plog$ and $s\in C^{\log}_{loc}(\rn)$.
\begin{list}{}{\labelwidth1.3em\leftmargin2em}
\item[{\upshape (i)\hfill}] If $p^+, q^+<\infty$, then the space $\Fspqpunkt$ is the collection of all $f\in \SSn$ such that 
\begin{equation*}
\|f \mid \Fspqpunkt\| := \Big\| \left(2^{js(\cdot)}(\varphi_j \widehat{f})^{\vee}\right)_{j \in \nat_0}\mid \ellpq\Big\|
\end{equation*}
is finite. 
\item[{\upshape (ii)\hfill}] The space $\Bspqpunkt$ is the collection of all $f\in \SSn$ such that 
\begin{equation*}
\|f \mid \Bspqpunkt\| := \Big\| \left(2^{js(\cdot)}(\varphi_j \widehat{f})^{\vee}\right)_{j \in \nat_0}\mid \ellqp\Big\|
\end{equation*}
is finite. 
\end{list}
\end{definition}

\begin{remark}
The independence of the resolution of unity in the definition of the spaces $\Fspqpunkt$ and $\Bspqpunkt$ can be justified by characterizations of the spaces, for example by local means (see \cite{Kem09} for Triebel-Lizorkin and \cite{KemV12} for Besov spaces), if $p,q\in\Plog$ and $s\in C^{\log}_{loc}(\rn)$. 
\end{remark}
\begin{remark}
These spaces include very well known spaces. In particular, if $\p=p$, $\q=q$ and $\s=s$ are constants, we get back to the classical Besov and Triebel-Lizorkin spaces $\Bspq$ and $\Fspq$. 
\end{remark}

The spaces $\Bspqpunkt$ and $\Fspqpunkt$ are isomorphic to sequence spaces. The underlying theorems are characterizations of the spaces above by atoms, wavelets and the $\varphi$-transform. The assertions can be found in \cite{DHR09}, \cite{Kem10} and \cite{Dri12}.\\
We do not repeat these characterizations here but we shall introduce some notation in order to state the sequence space characterizations. Let $\zn$ stand for the lattice of all points in $\rn$ with
integer-valued components,  $Q_{j,m}$ denotes a cube in $\rn$ with sides parallel to the axes of coordinates, centered at $2^{-j}m=(2^{-j}m_1,\dotsc,2^{-j}m_n)$ and with side length $2^{-j}$, where $m=(m_1,\dotsc,m_n)\in \zn$ and $j\in \nat_0$. If $Q$ is a cube in $\rn$ and $r>0$ then $rQ$ is the cube in $\rn$ concentric with $Q$ and with side length $r$ times the side length of $Q$. By $\chi_{E}$ we denote the characteristic function of the measurable set $E$. However, when $E$ is the cube  $Q_{j,m}$, the characteristic function of $Q_{j,m}$ is simply denoted by $\chi_{j,m}$.

\begin{definition}\label{defn:seq}
Let $p,q\in \Plog$ and $s\in C^{\log}_{loc}(\rn)$.

\begin{list}{}{\labelwidth1.3em\leftmargin2em}
\item[{\upshape (i)\hfill}] If $p^{+} < \infty$, then the sequence space $f_{\p,\q}^{\s}(\rn)$ consists of those complex-valued sequences $\lambda=(\lambda_{j,m})_{j\in\nat_0,m\in\zn}$ such that
\begin{align}
\|\lambda \mid f_{\p,\q}^{\s}(\rn)\|&:= \Big\Vert \Big(\sum_{m\in \zn}|\lambda_{j,m}|\,2^{j\s}\,\chi_{j,m} \Big)_{j\in\nat_0}\mid \ellpq \Big\Vert\notag
\end{align}
is finite. 
\item[{\upshape (ii)\hfill}] The sequence space $b_{\p,\q}^{\s}(\rn)$ consists of those complex-valued sequences $\lambda=(\lambda_{j,m})_{j\in\nat_0,m\in\zn}$ such that
\begin{align}
\|\lambda \mid b_{\p,\q}^{\s}(\rn)\|&:= \Big\Vert \Big(\sum_{m\in \zn}|\lambda_{j,m}|\,2^{j\s} \,\chi_{j,m} \Big)_{j\in\nat_0}\mid \ellqp \Big\Vert \notag
\end{align}
is finite.
\end{list}
\end{definition}

Now, we can state the connection between the function spaces and the corresponding sequence spaces. If we have a distribution $f\in\SSn$ then we can identify it with the corresponding sequence $\lambda(f)=\left(\lambda_{\nu,m}(f)\right)_{\nu\in\nat_0,m\in\zn}$ and vice versa. It depends on the underlying characterization (atoms, wavelets, $\varphi$-transform) how this connection is made, c.f. \cite{DHR09}, \cite{Dri12} and \cite{Kem10}.
\begin{theorem}\label{thm:characterization}
Let $p,q\in\Plog$ and $s\in C^{\log}_{loc}(\rn)$. 
\begin{enumerate}
\item We have
\begin{align*}
\norm{f}{\Bspqpunkt}\sim\norm{\lambda(f)}{b_{\p,\q}^{\s}(\rn)}.
\end{align*}
\item If $p^+,q^+<\infty$, then we have
\begin{align*}
\norm{f}{\Fspqpunkt}\sim\norm{\lambda(f)}{f_{\p,\q}^{\s}(\rn)}.
\end{align*}
\end{enumerate} 
The constants in both assertions are independent on $f\in\SSn$.
\end{theorem}

To prove our main results, we will make use of the Sobolev embedding for $\bspqpunkt$, proved in \cite{AH10}. The counterpart for $\fspqpunkt$ was proved by Vyb\'\i ral in \cite{Vyb09}.

\begin{theorem}\label{sobolev}
Let $p_0, p_1, q \in \Plog$ and $s_0, s_1\in C^{\log}_{loc}(\rn)$. Let $s_0(x) \geq s_1(x)$ and $p_0(x)\leq p_1(x)$ for all $x\in \rn$ with
$$s_0(x)-\frac{n}{p_0(x)} = s_1(x)-\frac{n}{p_1(x)}, \quad x\in \rn.$$ 
Then we have
        \beq
            b^{s_0(\cdot)}_{p_0(\cdot), \q}(\rn) \hookrightarrow b^{s_1(\cdot)}_{p_1(\cdot), \q}(\rn). \nonumber
        \eeq
%
\end{theorem}

\section{Franke embedding - constant exponents case}\label{sec:FrankeKonstant}
The Franke embedding \eqref{eq:FJ2} was shown in \cite{Fra86} using duality and interpolation.
An alternative proof avoiding interpolation was given in \cite{Vyb08}. The main tools used were
the technique of non-increasing rearrangements and duality.
Here, we give a new proof of the Franke embedding for the scale of spaces with constant exponents,
which still relies on non-increasing rearrangements, but we avoid using duality.

We start by introducing the concept of non-increasing rearrangement and some of its important properties.
We refer to \cite{BenSha} to an extensive treatment of this subject and also to the proofs of the lemmas given below.
\begin{definition}
 Let $\mu$ be the Lebesgue measure in $\rn$. If $h$ is a measurable function on $\rn$, we define the non-increasing rearrangement of $h$ through
 \beq
  h^*(t) = \sup \left\{\lambda >0 : \mu\{x\in \rn : |h(x)|>\lambda\}>t \right\}, \quad t\in (0, \infty). \nonumber
 \eeq
\end{definition}

\begin{lemma}\label{rea1}
If $0<p\leq \infty$, then
\beq
  \|h \mid L_p(\rn)\| = \|h^* \mid L_p(0, \infty)\|\nonumber
\eeq
for every measurable function $h$.
\end{lemma}

\begin{lemma}\label{rea2}
 Let $h_1$ and $h_2$ be two non-negative measurable functions on $\rn$. If $1\leq p \leq \infty$, then 
 \beq
  \|h_1+h_2 \mid L_p(\rn)\| \leq \|h_1^* + h_2^* \mid L_p(0, \infty)\|. \nonumber
\eeq
\end{lemma}

The main result of this section is an alternative proof of \eqref{eq:FJ2} for the sequence spaces of Besov and Triebel-Lizorkin type
with constant exponents.
\begin{theorem}\label{thm:FrankeKonstant} Let $-\infty<s_1<s_0<\infty$, $0<p_0<p_1<\infty$ and $0<q\leq \infty$ with
$$s_0-\frac{n}{p_0}=s_1-\frac{n}{p_1}.$$
Then
$$
b^{s_0}_{p_0,p_1}(\rn)\hookrightarrow f^{s_1}_{p_1,q}(\rn).
$$
\end{theorem}
\bpr
By lifting properties it is enough to show the Franke embedding with $s_1=0$, i.e.
\begin{equation}\label{eq:Fr1}
b^{n\bigl(\frac{1}{p_0}-\frac{1}{p_1}\bigr)}_{p_0,p_1}(\rn)\hookrightarrow f^{0}_{p_1,q}(\rn).
\end{equation}
Plugging in Definition \ref{defn:seq}, we get
\begin{align}\label{eq:Fr1:temp1}
\|\gamma|f^{0}_{p_1,q}(\rn)\|&=\Bigl\|\Bigl(\sum_{j=0}^\infty\sum_{m\in\zn}|\gamma_{j,m}|^q\chi_{j,m}(\cdot)\Bigr)^{\frac1q}|L_{p_1}(\rn)\Bigr\|
\end{align}
and
\begin{align}\label{eq:Fr1:temp2}
\|\gamma|b^{n\bigl(\frac{1}{p_0}-\frac{1}{p_1}\bigr)}_{p_0,p_1}(\rn)\|
= \Bigl(\sum_{j=0}^\infty 2^{-jn}\Bigl(\sum_{m\in\zn}|\gamma_{j,m}|^{p_0}\Bigr)^{\frac{p_1}{p_0}}\Bigr)^{\frac{1}{p_1}}.
\end{align}

To prove \eqref{eq:Fr1}, we will use $0<p_0<p_1<\infty$ and also $0<q\le \min(1,p_0)$, since it follows by elementary embeddings that we can take $q$ arbitrarily small.

We write
\begin{align*}
\|\gamma|f^{0}_{p_1,q}(\rn)\|&=\Bigl\|\Bigl(\sum_{j=0}^\infty\sum_{m\in\zn}|\gamma_{j,m}|^q\chi_{j,m}(\cdot)\Bigr)^{\frac1q}|L_{p_1}(\rn)\Bigr\|\\
&=\Bigl\|\sum_{j=0}^\infty\sum_{m\in\zn}|\gamma_{j,m}|^q\chi_{j,m}(\cdot)|L_{\frac{p_1}{q}}(\rn)\Bigr\|^{1/q}\\
&\le \Bigl\|\sum_{j=0}^\infty \Bigl(\sum_{m\in\zn}|\gamma_{j,m}|^q\chi_{j,m}\Bigr)^*(\cdot)|L_{\frac{p_1}{q}}(0,\infty)\Bigr\|^{1/q}\\
&= \Bigl\|\sum_{j=0}^\infty \sum_{l=0}^\infty(\gamma^*_{j,l})^q\chi^*_{j,l}(\cdot)|L_{\frac{p_1}{q}}(0,\infty)\Bigr\|^{1/q}.
\end{align*}
Here, we have used Lemma \ref{rea2} due to $\frac{p_1}{q}\ge 1$, $(\gamma^*_{j,l})_{l=0}^\infty$ is a non-increasing rearrangement of $(|\gamma_{j,m}|)_{m\in \zn}$
and $\chi^*_{j,l}$ is the characteristic function of $[2^{-jn}l,2^{-jn}(l+1)).$ Discretizing the last norm, we get by using the properties of $\chi^*_{j,l}(2^{-kn})$
\begin{align*}
\|\gamma|f^{0}_{p_1,q}(\rn)\|&\lesssim\Bigl(\sum_{k=-\infty}^{\infty}2^{-kn} \Bigl(\sum_{j=0}^\infty \sum_{l=0}^\infty(\gamma^*_{j,l})^q\chi^*_{j,l}(2^{-kn})\Bigr)^{\frac{p_1}{q}}\Bigr)^{\frac{1}{p_1}}\\
&=\Bigl(\sum_{k=-\infty}^{\infty}2^{-kn} \Bigl(\sum_{j=0}^\infty (\gamma^*_{j,\max(0,2^{(j-k)n})})^q\Bigr)^{\frac{p_1}{q}}\Bigr)^{\frac{1}{p_1}}\\
&\lesssim \Bigl[\sum_{k=-\infty}^{0}2^{-kn} \Bigl(\sum_{j=0}^\infty (\gamma^*_{j,2^{(j-k)n}})^q\Bigr)^{\frac{p_1}{q}}\Bigr]^{\frac{1}{p_1}}\\
&\qquad+\Bigl[\sum_{k=1}^\infty2^{-kn} \Bigl(\sum_{j=0}^{k-1} (\gamma^*_{j,0})^q\Bigr)^{\frac{p_1}{q}}\Bigr]^{\frac{1}{p_1}}\\
&\qquad+\Bigl[\sum_{k=1}^\infty2^{-kn} \Bigl(\sum_{j=k}^\infty (\gamma^*_{j,2^{(j-k)n}})^q\Bigr)^{\frac{p_1}{q}}\Bigr]^{\frac{1}{p_1}}=I+II+III.
\end{align*}
We estimate all the three terms separately.

The first term can be estimated in the following way. Starting with H{\"o}lder's inequality with $\beta>0$, we get
\begin{align*}
I&=\Bigl[\sum_{k=-\infty}^{0}2^{-kn} \Bigl(\sum_{j=0}^\infty (\gamma^*_{j,2^{(j-k)n}})^q\Bigr)^{\frac{p_1}{q}}\Bigr]^{\frac{1}{p_1}}\\
&\lesssim \Bigl[\sum_{k=-\infty}^{0} \Bigl(\sum_{j=0}^\infty 2^{-kn\frac{p_0}{p_1}}2^{jn\beta \frac{p_0}{q}}(\gamma^*_{j,2^{(j-k)n}})^{p_0}\Bigr)^{\frac{p_1}{p_0}}\Bigr]^{\frac{p_0}{p_1}\frac{1}{p_0}}.
\end{align*}
Now, after using the triangle inequality with $\frac{p_1}{p_0}>1$ and the embedding $\ell_{\frac{p_0}{p_1}}\hookrightarrow\ell_1$ in the sum over $k$, we substitute $l=j-k$ and get
\begin{align*}
I &\lesssim \Bigl[\sum_{j=0}^{\infty} \Bigl(\sum_{k=-\infty}^0 2^{-kn}2^{jn\beta \frac{p_1}{q}}(\gamma^*_{j,2^{(j-k)n}})^{p_1}\Bigr)^{\frac{p_0}{p_1}}\Bigr]^{\frac{1}{p_0}}\\
&\le \Bigl[\sum_{j=0}^{\infty} \sum_{k=-\infty}^0 2^{-kn \frac{p_0}{p_1}}2^{jn\beta \frac{p_0}{q}}(\gamma^*_{j,2^{(j-k)n}})^{p_0}\Bigr]^{\frac{1}{p_0}}\\
&= \Bigl[\sum_{j=0}^{\infty} 2^{jn\beta \frac{p_0}{q}}\sum_{l=j}^\infty 2^{(l-j)n \frac{p_0}{p_1}}(\gamma^*_{j,2^{ln}})^{p_0}\Bigr]^{\frac{1}{p_0}}\\
&= \Bigl[\sum_{j=0}^{\infty} 2^{jn(\beta \frac{p_0}{q}-\frac{p_0}{p_1})}\sum_{l=j}^\infty 2^{ln}2^{ln (\frac{p_0}{p_1}-1)}(\gamma^*_{j,2^{ln}})^{p_0}\Bigr]^{\frac{1}{p_0}}.
\end{align*}
Since $\frac{p_0}{p_1}-1 <0$, we have
\begin{align*}
I &\lesssim \Bigl[\sum_{j=0}^{\infty} 2^{jn(\beta \frac{p_0}{q}-\frac{p_0}{p_1})}2^{jn (\frac{p_0}{p_1}-1)}\sum_{l=j}^\infty 2^{ln}(\gamma^*_{j,2^{ln}})^{p_0}\Bigr]^{\frac{1}{p_0}}\\
&= \Bigl[\sum_{j=0}^{\infty} 2^{-jn\frac{p_0}{p_1}}2^{jn(\beta\frac{p_0}{q}+\frac{p_0}{p_1}-1)}\sum_{l=j}^\infty 2^{ln}(\gamma^*_{j,2^{ln}})^{p_0}\Bigr]^{\frac{1}{p_0}}.
\end{align*}
We finish this estimate choosing $\beta$ with $\beta\frac{p_0}{q}+\frac{p_0}{p_1}-1<0$, i.e., $0<\beta<q(\frac{1}{p_0}-\frac{1}{p_1})$ and get with \eqref{eq:Fr1:temp2}
\begin{align*}
I&\lesssim \Bigl[\sum_{j=0}^{\infty} 2^{-jn}\Bigl(\sum_{l=j}^\infty 2^{ln}(\gamma^*_{j,2^{ln}})^{p_0}\Bigr)^{\frac{p_1}{p_0}}\Bigr]^{\frac{1}{p_1}}
\lesssim\|\gamma|b^{n\bigl(\frac{1}{p_0}-\frac{1}{p_1}\bigr)}_{p_0,p_1}(\rn)\|.
\end{align*}

To estimate the second term, we start by using H{\"o}lder's inequality with $p_1/q>1$
\begin{align*}
II^{p_1}
&\lesssim \sum_{k=1}^\infty 2^{-kn}2^{kn\beta \frac{p_1}{q}}\sum_{j=0}^{k-1}2^{-jn\beta \frac{p_1}{q}}(\gamma_{j,0}^*)^{{p_1}}\\
&=\sum_{j=0}^\infty 2^{-jn\beta\frac{p_1}{q}}(\gamma^*_{j,0})^{p_1} \sum_{k=j+1}^\infty 2^{kn(\beta \frac{p_1}{q}-1)}.
\end{align*}
Now, choosing $\beta>0$ with $\beta \frac{p_1}{q}-1<0$, i.e., $0<\beta<\frac{q}{p_1}$, we get again by \eqref{eq:Fr1:temp2}
\begin{align*}
II^{p_1}&\lesssim\sum_{j=0}^\infty 2^{-jn\beta\frac{p_1}{q}}(\gamma^*_{j,0})^{p_1}2^{jn\beta\frac{p_1}{q}}2^{-jn}\\
&=\sum_{j=0}^\infty (\gamma^*_{j,0})^{p_1}2^{-jn}\lesssim\|\gamma|b^{n\bigl(\frac{1}{p_0}-\frac{1}{p_1}\bigr)}_{p_0,p_1}(\rn)\|^{p_1}.
\end{align*}

Finally, we estimate the third term starting again by using a parameter $\delta>0$ and H{\"o}lder's inequality
\begin{align*}
III^{p_1} 
&=\sum_{k=1}^\infty2^{-kn} \Bigl(\sum_{j=k}^\infty 2^{jn\delta}2^{-jn\delta}(\gamma^*_{j,2^{(j-k)n}})^q\Bigr)^{\frac{p_1}{q}}\\
&\lesssim \sum_{k=1}^\infty2^{-kn} 2^{-kn\delta \frac{p_1}{q}}\sum_{j=k}^\infty 2^{jn\delta \frac{p_1}{q}}(\gamma^*_{j,2^{(j-k)n}})^{p_1}\\
&= 
\sum_{j=1}^\infty 2^{jn\delta \frac{p_1}{q}}\sum_{k=1}^j 2^{-kn(1+\delta \frac{p_1}{q})}(\gamma^*_{j,2^{(j-k)n}})^{p_1}.
\end{align*}
We substitute $l=k-j$ and obtain
\begin{align*}
III^{p_1} 
&\lesssim \sum_{j=1}^\infty 2^{jn\delta \frac{p_1}{q}}\sum_{l=0}^{j-1} 2^{(l-j)n(1+\delta \frac{p_1}{q})}(\gamma^*_{j,2^{ln}})^{p_1}\\
&=\sum_{j=1}^\infty 2^{-jn}\Bigl[\sum_{l=0}^{j-1} 2^{ln(\frac{q}{p_1}+\delta )\frac{p_1}{q}}(\gamma^*_{j,2^{ln}})^{p_1}\Bigr]\\
&\le\sum_{j=1}^\infty 2^{-jn}\Bigl[\sum_{l=0}^{j-1} 2^{ln(\frac{q}{p_1}+\delta )\frac{p_0}{q}}(\gamma^*_{j,2^{ln}})^{{p_0}}\Bigr]^{\frac{p_1}{p_0}},
\end{align*}
where the last step comes from the elementary embedding $\ell_{p_0}\hookrightarrow \ell_{p_1}$. Choosing now $\delta>0$ such that $(\frac{q}{p_1}+\delta)\frac{p_0}{q}=1$, i.e. $\delta=\frac{q}{p_0}-\frac{q}{p_1}>0$, we get
\begin{align*}
III^{p_1} &\lesssim \sum_{j=1}^\infty 2^{-jn}\Bigl[\sum_{l=0}^{j-1} 2^{ln}(\gamma^*_{j,2^{ln}})^{p_0}\Bigr]^{\frac{p_1}{p_0}}
\lesssim\|\gamma|b^{n\bigl(\frac{1}{p_0}-\frac{1}{p_1}\bigr)}_{p_0,p_1}(\rn)\|^{p_1},
\end{align*}
which concludes the proof. 
\epr

\section{Franke-Jawerth embeddings - variable exponent case}\label{sec:FrankeJawerth}
We now return to the scale of Besov and Triebel-Lizorkin spaces with variable exponents. First, we state the result in the form of embeddings of sequence spaces
under the assumptions we really need in the proof. After that, we combine those with the conditions required in Theorem \ref{thm:characterization},
and we present the results on function spaces in the form of a corollary. 
\subsection{Jawerth embedding}
\begin{theorem} Let $p_0, p_1, q \in \Plog$ with $p_0^+<\infty$ and $s_0, s_1\in C^{\log}_{loc}(\rn)$.
Let 
$\displaystyle \inf_{x\in\rn}(s_0(x)-s_1(x))>0$
with
\begin{equation}\label{eq:Jaw:4}
s_0(x)-\frac{n}{p_0(x)} = s_1(x)-\frac{n}{p_1(x)}, \quad x\in \rn.
\end{equation}
Then
\begin{equation}\label{eq:Jaw:1}
f^{s_0(\cdot)}_{p_0(\cdot),\q}(\rn)\hookrightarrow b^{s_1(\cdot)}_{p_1(\cdot),p_0(\cdot)}(\rn).
\end{equation}
\end{theorem}
\bpr
Let us put
$$
\varepsilon'=\inf_{x\in\rn}(s_0(x)-s_1(x))=\inf_{x\in\rn}\Bigl(\frac{n}{p_0(x)}-\frac{n}{p_1(x)}\Bigr)>0.
$$
Then
\begin{align*}
\frac{p_1(x)}{p_0(x)}-1=p_1(x)\Bigl(\frac{1}{p_0(x)}-\frac{1}{p_1(x)}\Bigr)\ge p_1(x)\frac{\varepsilon'}{n}\ge \frac{p_1^-\varepsilon'}{n}.
\end{align*}
Putting $\varepsilon = \frac{p_1^-\varepsilon'}{2n}>0$ we get for every $x\in\rn$
\begin{equation}\label{eq:Jaw:3}
p_0(x)<(1+\varepsilon)p_0(x)< p_1(x).
\end{equation}
The proof of \eqref{eq:Jaw:1} will be the result of the following chain of embeddings
\begin{align}
\notag f^{s_0(\cdot)}_{p_0(\cdot),\q}(\rn)&\hookrightarrow f^{s_0(\cdot)}_{p_0(\cdot),\infty}(\rn)
\hookrightarrow b^{s_0(\cdot)-\frac{n}{p_0(\cdot)}+\frac{n}{(1+\varepsilon)p_0(\cdot)}}_{(1+\varepsilon)p_0(\cdot),p_0(\cdot)}(\rn)\\
\label{eq:Jaw:2}&=b^{s_1(\cdot)-\frac{n}{p_1(\cdot)}+\frac{n}{(1+\varepsilon)p_0(\cdot)}}_{(1+\varepsilon)p_0(\cdot),p_0(\cdot)}(\rn)
\hookrightarrow b^{s_1(\cdot)}_{p_1(\cdot),p_0(\cdot)}(\rn).
\end{align}
The first embedding in \eqref{eq:Jaw:2} is an elementary statement about the monotonicity of $f$-spaces in the summability index $q$
and the last embedding follows from Theorem \ref{sobolev} and \eqref{eq:Jaw:3}. The identity in \eqref{eq:Jaw:2}
is a simple consequence of \eqref{eq:Jaw:4}. Hence, it remains to prove the second embedding in \eqref{eq:Jaw:2},
which is actually a special case of \eqref{eq:Jaw:1} with $q=\infty$ and $p_1(x)=(1+\varepsilon)p_0(x).$
Finally, by the lifting property of Besov and Triebel-Lizorkin spaces with variable exponents, we may consider only the case
when the smoothness exponent of the target space is zero.

The proof of \eqref{eq:Jaw:1} will therefore follow from
\begin{equation}\label{eq:want}
f^{\frac{n}{p_0(\cdot)}\cdot\frac{\varepsilon}{1+\varepsilon}}_{p_0(\cdot),\infty}(\rn)\hookrightarrow b^0_{(1+\varepsilon)p_0(\cdot),p_0(\cdot)}(\rn)
\end{equation}
or, equivalently,
\begin{equation} \label{norm:want1}
\|\gamma \mid b^{0}_{(1+\varepsilon)p_0(\cdot),p_0(\cdot)} (\rn)\| \leq C \, \|\gamma \mid f^{\frac{n}{p_0(\cdot)}\cdot\frac{\varepsilon}{1+\varepsilon}}_{p_0(\cdot),\infty}(\rn) \|,
\end{equation}
for some constant $C>0$ and $\gamma=(\gamma_{j,m})_{j,m}$, for $j\in\nat_0, m\in\zn$.

Let us put
$$
h(x)=\sup_{j, m}2^{\frac{jn}{p_0(x)}\cdot\frac{\varepsilon}{1+\varepsilon}}|\gamma_{j,m}|\chi_{j,m}(x),\quad x\in\rn.
$$
Then for every $x\in Q_{j,m}$ we have $2^{\frac{jn}{p_0(x)}\cdot\frac{\varepsilon}{1+\varepsilon}}|\gamma_{j,m}|\le h(x)$ and
$$
|\gamma_{j,m}|\le \inf_{y\in Q_{j,m}}2^{\frac{-jn}{p_0(y)}\cdot\frac{\varepsilon}{1+\varepsilon}}h(y).
$$
Using this notation,
\begin{align*}
\|\gamma|f^{\frac{n}{p_0(x)}\cdot\frac{\varepsilon}{1+\varepsilon}}_{p_0(\cdot),\infty}(\rn)\|=\|h \mid L_{p_0(\cdot)}(\rn)\|
\end{align*}
and \eqref{norm:want1} reads as
\begin{equation} 
\|\gamma \mid b^{0}_{(1+\varepsilon)p_0(\cdot),p_0(\cdot)} (\rn)\| \leq C \, \|h \mid L_{p_0(\cdot)} (\rn)\|. \nonumber
\end{equation}

We assume that $\|h \mid L_{p_0(\cdot)} \| \leq 1$, i.e.
\begin{equation}\label{eq:h}
    \int_{\rn} h(x)^{p_0(x)} \, dx \leq 1,
\end{equation}
and need to prove that
\begin{align*}
\|\gamma \mid b^{0}_{(1+\varepsilon)p_0(\cdot),p_0(\cdot)}(\rn) \|& \leq C,
\intertext{which we will show by}
\varrho_{\ell_{p_0(\cdot)}(L_{(1+\varepsilon)p_0(\cdot)})}\left(\sum_{m \in \zn} |\gamma_{j,m}| \chi_{j,m}\right)&\leq C.
\end{align*}
We have
\begin{align*}
&\varrho_{\ell_{p_0(\cdot)}(L_{(1+\varepsilon)p_0(\cdot)})}\left(\sum_{m \in \zn} |\gamma_{j,m}| \chi_{j,m}\right)
= \sum_{j=0}^{\infty} \Bigl\| \sum_{m \in \zn} |\gamma_{j,m}|^{p_0(\cdot)} \chi_{j,m} \mid L_{\frac{(1+\varepsilon)p_0(\cdot)}{p_0(\cdot)}}(\rn)\Bigr\|\\
& \leq \sum_{j=0}^{\infty} \Bigl\| \sum_{m \in \zn}  \Bigl(\inf_{y\in Q_{j,m}}2^{\frac{-jn}{p_0(y)}\cdot\frac{\varepsilon}{1+\varepsilon}}h(y)\Bigr)^{p_0(\cdot)}
\chi_{j,m} \mid L_{1+\varepsilon}(\rn)\Bigr\|\\
&=\sum_{j=0}^{\infty} \Bigl\{\sum_{m \in \zn} \int_{Q_{j,m}} \Bigl(\inf_{y\in Q_{j,m}}2^{\frac{-jn}{p_0(y)}\cdot\frac{\varepsilon}{1+\varepsilon}}h(y)\Bigr)^{(1+\varepsilon)p_0(x)}dx\Bigr\}^{\frac{1}{1+\varepsilon}}\\
&=\sum_{j=0}^{\infty} 2^{-jn\frac{\varepsilon}{1+\varepsilon}}\Bigl\{\sum_{m \in \zn} \int_{Q_{j,m}}
\Bigl(\inf_{y\in Q_{j,m}}2^{-jn\varepsilon\bigl(\frac{p_0(x)}{p_0(y)}-1\bigr)}h(y)^{(1+\varepsilon)p_0(x)}\Bigr)dx\Bigr\}^{\frac{1}{1+\varepsilon}}.
\end{align*}
We use the regularity of $p_0$ to obtain for $x,y\in Q_{j,m}$
\begin{align*}
2^{-jn\bigl(\frac{p_0(x)}{p_0(y)}-1\bigr)}&
=\Bigl[2^{j\bigl(\frac{1}{p_0(x)}-\frac{1}{p_0(y)}\bigr)}\Bigr]^{np_0(x)}\le 2^{np_0(x)c_{\rm log}(1/p_0)}\le c'.
\end{align*}
So, it is enough to prove
\beq \label{eq:pos5}
\sum_{j=0}^\infty 2^{-jn\frac{\varepsilon}{1+\varepsilon}}\Bigl\{\sum_{m\in \zn}\int_{Q_{j,m}}h_{j,m}^{(1+\varepsilon)p_0(x)}dx\Bigr\}^{\frac{1}{1+\varepsilon}} \leq C,
\eeq
where we denoted
$$
h_{j,m}=\inf_{y\in Q_{j,m}}h(y).
$$

We split the left-hand side of \eqref{eq:pos5} into
\begin{align*}
\sum_{j=0}^\infty &2^{-jn\frac{\varepsilon}{1+\varepsilon}}\Bigl\{\sum_{m\in \zn}\int_{Q_{j,m}}h_{j,m}^{(1+\varepsilon)p_0(x)}dx\Bigr\}^{\frac{1}{1+\varepsilon}}\\
&\le\sum_{j=0}^\infty 2^{-jn\frac{\varepsilon}{1+\varepsilon}}\Bigl\{\sum_{\{m:h_{j,m}\le 1\}}\int_{Q_{j,m}}h_{j,m}^{(1+\varepsilon)p_0(x)}dx\Bigr\}^{\frac{1}{1+\varepsilon}} \\
&\quad +\sum_{j=0}^\infty 2^{-jn\frac{\varepsilon}{1+\varepsilon}}\Bigl\{\sum_{\{m:h_{j,m}>1\}}\int_{Q_{j,m}}h_{j,m}^{(1+\varepsilon)p_0(x)}dx\Bigr\}^{\frac{1}{1+\varepsilon}}=I+II.
\end{align*}

The first term can be estimated by \eqref{eq:h}
\begin{align*}
I
&\le \sum_{j=0}^\infty 2^{-jn\frac{\varepsilon}{1+\varepsilon}}\Bigl\{\sum_{\{m:h_{j,m}\le 1\}}\int_{Q_{j,m}}h_{j,m}^{p_0(x)}dx\Bigr\}^{\frac{1}{1+\varepsilon}}\\
&\le \sum_{j=0}^\infty 2^{-jn\frac{\varepsilon}{1+\varepsilon}}\Bigl\{\int_{\rn}h(x)^{p_0(x)}dx\Bigr\}^{\frac{1}{1+\varepsilon}}\le c.
\end{align*}

To estimate $II$, we make first a couple of observations. Let $h_{j,m}\ge 1$. Then
\begin{equation*}
1\ge \int_{\rn}h(x)^{p_0(x)}dx\ge \int_{Q_{j,m}}h(x)^{p_0(x)}dx\ge \int_{Q_{j,m}} h_{j,m}^{p_0(x)}dx
\ge |Q_{j,m}|h_{j,m}^{p_0^-}=2^{-jn}h_{j,m}^{p_0^-}
\end{equation*}
and we get
$$
1\le h_{j,m}\le 2^{jn/p_0^-}.
$$
Hence there is an $0\le \alpha\le 1$ such that $h_{j,m}=2^{\alpha jn/p_0^-}.$ Since $p_0^+<\infty$, we have that $1/p_0\in C^{\log}_{loc}(\rn)$ implies $p_0\in C^{\log}_{loc}(\rn)$ and we can use its regularity. Hence, there is a constant $c>1$ such that for any $x,y\in Q_{j,m}$ it holds $c^{-1}\le 2^{j(p_0(x)-p_0(y))}\le c$
and therefore
\begin{equation}\label{eq:help1}
c^{-\alpha n/p_0^-}\le 2^{\frac{\alpha jn}{p_0^-}(p_0(x)-p_0(y))}=h_{j,m}^{(p_0(x)-p_0(y))}\le c^{\alpha n/p_0^-}\le c^{n/p_0^-}.
\end{equation}
If we also denote $p_{j,m}=\inf_{y\in Q_{j,m}}p_0(y)$, we obtain
\begin{equation*}
h_{j,m}^{p_0(x)}=h_{j,m}^{p_0(x)-p_{j,m}}\cdot h_{j,m}^{p_{j,m}}\le C h_{j,m}^{p_{j,m}}\le C \inf_{y\in Q_{j,m}}\bigl(h(y)^{p_0(y)}\bigr).
\end{equation*}
The last fact we shall use is that
\begin{equation}\label{eq:this}
\sum_{j=0}^\infty 2^{-jn\frac{\varepsilon}{1+\varepsilon}}\Bigl\{\sum_{m\in \zn}\int_{Q_{j,m}}(\inf_{y\in Q_{j,m}}\varphi(y))^{1+\varepsilon}dx\Bigr\}^{\frac{1}{1+\varepsilon}}\le c\|\varphi \mid L_1(\rn)\|
\end{equation}
for each $\varphi\in L_1(\rn)$. The proof follows easily using the technique of non-increasing rearrangement.
\begin{align*}
\sum_{j=0}^\infty &2^{-jn\frac{\varepsilon}{1+\varepsilon}}\Bigl\{\sum_{m\in \zn}\int_{Q_{j,m}}(\inf_{y\in Q_{j,m}}\varphi(y))^{1+\varepsilon}dx\Bigr\}^{\frac{1}{1+\varepsilon}}\\
&\le \sum_{j=0}^\infty 2^{-jn\frac{\varepsilon}{1+\varepsilon}}\Bigl\{\sum_{l=1}^\infty 2^{-jn}\varphi^*(l2^{-jn})^{1+\varepsilon}\Bigr\}^{\frac{1}{1+\varepsilon}}\\
&\lesssim \sum_{j=0}^\infty 2^{-jn\frac{\varepsilon}{1+\varepsilon}}\Bigl\{\sum_{k=0}^\infty 2^{(k-j)n}\varphi^*(2^{(k-j)n})^{1+\varepsilon}\Bigr\}^{\frac{1}{1+\varepsilon}}\\
&\le \sum_{j=0}^\infty 2^{-jn\frac{\varepsilon}{1+\varepsilon}}\sum_{k=0}^\infty 2^{(k-j)\frac{n}{1+\varepsilon}}\varphi^*(2^{(k-j)n})\\
&= \sum_{l=-\infty}^\infty \varphi^*(2^{-ln})2^{-\frac{ln}{1+\varepsilon}}\sum_{j=l}^\infty 2^{-jn\frac{\varepsilon}{1+\varepsilon}}\lesssim
\sum_{l=-\infty}^\infty \varphi^*(2^{-ln})2^{-ln}\sim\|\varphi \mid L_1(\rn)\|.
\end{align*}
We apply \eqref{eq:this} with $\varphi(y)=h(y)^{p_0(y)}$ to estimate II
\begin{align*}
II&=\sum_{j=0}^\infty 2^{-jn\frac{\varepsilon}{1+\varepsilon}}\Bigl\{\sum_{\{m:h_{j,m}>1\}}\int_{Q_{j,m}}h_{j,m}^{(1+\varepsilon)p_0(x)}dx\Bigr\}^{\frac{1}{1+\varepsilon}}\\
&\lesssim\sum_{j=0}^\infty 2^{-jn\frac{\varepsilon}{1+\varepsilon}}\Bigl\{\sum_{\{m:h_{j,m}>1\}}\int_{Q_{j,m}}(\inf_{y\in Q_{j,m}}h(y)^{p_0(y)})^{1+\varepsilon}dx\Bigr\}^{\frac{1}{1+\varepsilon}}\\
&\lesssim\sum_{j=0}^\infty 2^{-jn\frac{\varepsilon}{1+\varepsilon}}\Bigl\{\sum_{m\in\zn}\int_{Q_{j,m}}(\inf_{y\in Q_{j,m}}h(y)^{p_0(y)})^{1+\varepsilon}dx\Bigr\}^{\frac{1}{1+\varepsilon}}\\
&\lesssim\int_{\rn}h(y)^{p_0(y)}dy\le C
\end{align*}
and finish the proof.
\epr
Using the correspondence of sequence and function spaces from Theorem \ref{thm:characterization} we obtain the Jawerth embedding for the variable function spaces.
\begin{corollary} Let $p_0, p_1, q \in \Plog$ with $p_0^+,q^+<\infty$ and $s_0, s_1\in C^{\log}_{loc}(\rn)$. Let 
$\inf_{x\in\rn}(s_0(x)-s_1(x))>0$ 
with $$s_0(x)-\frac{n}{p_0(x)} = s_1(x)-\frac{n}{p_1(x)}, \quad x\in \rn.$$
Then
$$
F^{s_0(\cdot)}_{p_0(\cdot),\q}(\rn)\hookrightarrow B^{s_1(\cdot)}_{p_1(\cdot),p_0(\cdot)}(\rn).
$$
\end{corollary}

\subsection{Franke embedding}
In this section we prove the Franke embedding for function spaces with variable exponents.
We have to avoid duality arguments in the variable exponent setting and therefore reduce the proof to the constant exponent case and apply Theorem \ref{thm:FrankeKonstant}, which we have shown in the previous section. 
\begin{theorem} Let $p_0, p_1, q \in \Plog$ with $p_1^+<\infty$ and $s_0, s_1\in C^{\log}_{loc}(\rn)$. Let 
$\inf_{x\in\rn}(s_0(x)-s_1(x))>0$ 
with
$$s_0(x)-\frac{n}{p_0(x)} = s_1(x)-\frac{n}{p_1(x)}, \quad x\in \rn.$$
Then
$$
b^{s_0(\cdot)}_{p_0(\cdot),p_1(\cdot)}(\rn)\hookrightarrow f^{s_1(\cdot)}_{p_1(\cdot),\q}(\rn).
$$
\end{theorem}
\bpr By the lifting property we may suppose again $s_1=0$, and the elementary embeddings between $f$-spaces allow one to set
\begin{equation}
 q(x)= \frac{1}{r}p_1(x), \quad \mbox{with } r>1 \mbox{ chosen big enough.} \nonumber
\end{equation}
So it suffices to prove the embedding
\begin{equation}\label{eq:pos6}
b^{n(\frac{1}{p_0(\cdot)}-\frac{1}{p_1(\cdot)})}_{p_0(\cdot),p_1(\cdot)}(\rn)\hookrightarrow f^{0}_{p_1(\cdot),\frac{1}{r}p_1(\cdot)}(\rn).
\end{equation}
Similarly to \eqref{eq:Jaw:3}, $p_0(x)<p_1(x)$ are again well separated and we may find $\varepsilon >0$ with
\begin{equation}
p_1(x)>(1-\varepsilon)p_1(x)>p_0(x). \nonumber
\end{equation}
By the Sobolev embedding from Theorem \ref{sobolev} we obtain
\begin{equation}
b^{n(\frac{1}{p_0(\cdot)}-\frac{1}{p_1(\cdot)})}_{p_0(\cdot),p_1(\cdot)}(\rn)\hookrightarrow b^{\frac{n}{p_1(\cdot)}\frac{\varepsilon}{1-\varepsilon}}_{(1-\varepsilon)p_1(\cdot),p_1(\cdot)}(\rn). \nonumber
\end{equation}
Hence, instead of \eqref{eq:pos6} we show
\begin{equation}\label{eq:pos7}
b^{\frac{n}{p_1(\cdot)}\frac{\varepsilon}{1-\varepsilon}}_{(1-\varepsilon)p_1(\cdot),p_1(\cdot)}(\rn)\hookrightarrow f^0_{p_1(\cdot),\frac{1}{r}p_1(\cdot)}(\rn).
\end{equation}

We assume that $\| \gamma \mid b^{\frac{n}{p_1(\cdot)}\frac{\varepsilon}{1-\varepsilon}}_{(1-\varepsilon)p_1(\cdot),p_1(\cdot)}(\rn)\| \leq 1$, which is equivalent to
\begin{equation}
 \sum_{j=0}^{\infty} \Bigl\| \sum_{m\in \zn} |\gamma_{j,m}|^{p_1(\cdot)}2^{jn\frac{\varepsilon}{1-\varepsilon}} \chi_{j,m} \mid L_{1-\varepsilon}(\rn)\Bigr\| \leq 1 \nonumber
\end{equation}
or even
\begin{equation}\label{eq:pos7a}
 \sum_{j=0}^{\infty}2^{jn\frac{\varepsilon}{1-\varepsilon}} \left( \int_{\rn} \sum_{m\in \zn} |\gamma_{j,m}|^{(1-\varepsilon)p_1(x)} \chi_{j,m}(x) dx\right)^{\frac{1}{1-\varepsilon}} \leq 1.
\end{equation}
This implies that, for every $j \in \nat_0$, 
we have the inequality
\begin{equation}\label{eq:pos8}
\sum_{m\in\zn} \int_{Q_{j,m}}  |\gamma_{j,m}|^{(1-\varepsilon)p_1(x)}  dx  \leq2^{-jn\varepsilon}.
\end{equation}

Our aim is to prove that $\|\gamma \mid f^0_{p_1(\cdot),\frac{1}{r}p_1(\cdot)}(\rn)\| \leq C$, which is equivalent to prove the following
\begin{equation}\label{eq:pos9}
 \Bigl\| \sum_{j=0}^{\infty} \sum_{m\in \zn} |\gamma_{j,m}|^{\frac{p_1(\cdot)}{r}} \chi_{j,m} \mid L_r (\rn)\Bigr\| \leq C.
\end{equation}

We have
\begin{align*}
\Bigl\| &\sum_{j=0}^{\infty}  \sum_{m\in \zn} |\gamma_{j,m}|^{\frac{p_1(\cdot)}{r}} \chi_{j,m} \mid L_r(\rn) \Bigr\| \\
&\leq \Bigl\|\sum_{j=0}^{\infty} \sum_{\{m:|\gamma_{j,m}|\le 1\}} |\gamma_{j,m}|^{\frac{p_1(\cdot)}{r}}  \chi_{j,m} \mid L_r(\rn) \Bigr\|
+ \Bigl\|\sum_{j=0}^{\infty}\sum_{\{m:|\gamma_{j,m}|>1\}} |\gamma_{j,m}|^{\frac{p_1(\cdot)}{r}}  \chi_{j,m} \mid L_r(\rn) \Bigr\| \\
&=I+II.
\end{align*}
The estimate of $I$ follows by \eqref{eq:pos8}
\begin{align*}
 I &\leq \sum_{j=0}^\infty \Bigl(\sum_{\{m:|\gamma_{j,m}|\le 1\}}\int_{Q_{j,m}}|\gamma_{j,m}|^{p_1(x)}dx\Bigr)^{\frac{1}{r}} \leq \sum_{j=0}^\infty \Bigl(\sum_{\{m:|\gamma_{j,m}|\le 1\}}\int_{Q_{j,m}}|\gamma_{j,m}|^{(1-\varepsilon)p_1(x)}dx\Bigr)^{\frac{1}{r}} \\
 & \leq \sum_{j=0}^\infty 2^{-jn\varepsilon \frac1r} \leq c.
\end{align*}
To estimate $II$, we observe that \eqref{eq:pos8} implies for every $j\in\nat_0$ and $m\in\zn$ with $|\gamma_{j,m}|> 1$
\begin{equation*}
 2^{-jn\varepsilon} \geq \int_{Q_{j,m}}  |\gamma_{j,m}|^{(1-\varepsilon)p_1(x)} dx \geq 2^{-jn} |\gamma_{j,m}|^{(1-\varepsilon)p_1^-}
\end{equation*}
and therefore
\begin{equation*}
 1\leq|\gamma_{j,m}| \leq 2^{jn \frac{1}{p_1^-}}.
\end{equation*}
Similarly to \eqref{eq:help1}, for every $x\in Q_{j,m}$ we have $|\gamma_{j,m}|^{p_1(x)}\sim |\gamma_{j,m}|^{p_{j,m}}$,
where $p_{j,m}$ is the value of $p_1$ in the middle of $Q_{j,m}$. Denoting $\alpha_{j,m}=|\gamma_{j,m}|^{p_{j,m}}$,
we get that \eqref{eq:pos7a} implies
\begin{equation}\label{eq:last6}
\sum_{j=0}^\infty 2^{jn\frac{\varepsilon}{1-\varepsilon}}\Bigl(2^{-jn}\sum_{\{m:|\gamma_{j,m}|> 1\}}|\gamma_{j,m}|^{(1-\varepsilon)p_{j,m}}\Bigr)^{\frac{1}{1-\varepsilon}}\le c,
\end{equation}
or, equivalently,
\begin{equation}\label{eq:last7}
\sum_{j=0}^\infty 2^{-jn}\Bigl(\sum_{\{m:|\alpha_{j,m}|> 1\}}|\alpha_{j,m}|^{1-\varepsilon}\Bigr)^{\frac{1}{1-\varepsilon}}\le c.
\end{equation}
To estimate the second term II it is therefore sufficient to show that
\begin{equation}\label{eq:last8}
\Bigl\|\sum_{j=0}^\infty\sum_{\{m:|\gamma_{j,m}|> 1\}}|\gamma_{j,m}|^{p_{j,m}/r}\chi_{j,m}|L_r(\rn)\Bigr\|\le C,
\end{equation}
which in turn is equivalent to
\begin{equation}\label{eq:last9}
\Bigl\|\sum_{j=0}^\infty\sum_{\{m:|\alpha_{j,m}|> 1\}}|\alpha_{j,m}|^{1/r}\chi_{j,m}|L_r(\rn)\Bigr\|\le C.
\end{equation}
To obtain \eqref{eq:last9} from \eqref{eq:last7}, we employ the constant-index case, i.e. Theorem \ref{thm:FrankeKonstant}.
Indeed, we put $\beta_{j,m}=|\alpha_{j,m}|^{1/r}$ if $|\alpha_{j,m}|>1$ and zero otherwise and get
\begin{align*}
\Bigl\|\sum_{j=0}^\infty&\sum_{\{m:|\alpha_{j,m}|> 1\}}|\alpha_{j,m}|^{1/r}\chi_{j,m}|L_r(\rn)\Bigr\|=
\Bigl\|\sum_{j=0}^\infty\sum_{m\in\zn}|\beta_{j,m}|\chi_{j,m}|L_r(\rn)\Bigr\|\\
&=\|\beta\mid f^{0}_{r,1}\|\lesssim\|\beta\mid b^{\frac{n\varepsilon}{r(1-\varepsilon)}}_{(1-\varepsilon)r,r}\|
=\Bigl\{\sum_{j=0}^\infty 2^{\frac{jn\varepsilon}{(1-\varepsilon)r}r}\Bigl(\sum_{m\in\zn}2^{-jn}\beta_{j,m}^{(1-\varepsilon)r}\Bigr)^{\frac{r}{(1-\varepsilon)r}}\Bigr\}^{1/r}\\
&=\Bigl\{\sum_{j=0}^\infty 2^{-jn}\Bigl(\sum_{m:|\alpha_{j,m}|>1}|\alpha_{j,m}|^{1-\varepsilon}\Bigr)^{\frac{1}{1-\varepsilon}}\Bigr\}^{1/r}\le c^{1/r}.
\end{align*}
\epr
\begin{corollary} Let $p_0, p_1, q \in \Plog$ with $p_1^+, q^+<\infty$ and $s_0, s_1\in C^{\log}_{loc}(\rn)$. Let 
$\inf_{x\in\rn}(s_0(x)-s_1(x))>0$ 
with
$$s_0(x)-\frac{n}{p_0(x)} = s_1(x)-\frac{n}{p_1(x)}, \quad x\in \rn.$$
Then
$$
B^{s_0(\cdot)}_{p_0(\cdot),p_1(\cdot)}(\rn)\hookrightarrow F^{s_1(\cdot)}_{p_1(\cdot),\q}(\rn).
$$
\end{corollary}

\section{Jawerth and Franke embedding in 2-microlocal spaces}\label{sec:2ml}
The definition of Besov and Triebel-Lizorkin spaces of variable smoothness and integrability is a special case of the so-called 2-microlocal spaces of variable integrability. 
As all the proofs for spaces of variable smoothness do also serve for 2-microlocal spaces, we devote this chapter to present these results. We start by the definition of the spaces, which is based on the dyadic decomposition of unity as presented before combined with the concept of admissible weight sequences.

\begin{definition}\label{def-ad-weight} Let $\alpha\geq 0$ and $\alpha_1,\alpha_2\in\real$ with $\alpha_1\leq \alpha_2$. A
sequence of non-negative measurable functions in $\rn$ $\bm{w}=(w_j)_{j\in\nat_0}$ belongs to the class $\mathcal{W}^{\alpha}_{\alpha_1,\alpha_2}(\rn)$ if the following conditions are satisfied:
  \begin{list}{}{\labelwidth1.7em\leftmargin2.3em}
    \item[{\hfill (i)\hfill}] There exists a constant $c>0$ such that
      $$
	0<w_j(x)\leq c\,w_j(y)\,(1+2^j|x-y|)^{\alpha}\quad \text{for all} \;\, j\in\nat_0 \;\; \text{and all} \;\, x,y\in\rn.
      $$
    \item[{\hfill (ii)\hfill}] For all $j\in\nat_0$ it holds
      $$
	2^{\alpha_1}\,w_j(x)\leq w_{j+1}(x)\leq 2^{\alpha_2}\,w_j(x) \quad \text{for all}\;\, x \in \rn.
      $$
  \end{list}
Such a system $(w_j)_{j\in\nat_0}\in\mathcal{W}^{\alpha}_{\alpha_1,\alpha_2}(\rn)$ is called admissible weight sequence.\\
\end{definition}

Properties of admissible weights may be found in \cite[Remark~2.4]{Kem08}. Finally, here is the definition of the spaces under consideration.

\begin{definition}  \label{def-Wspaces}
Let $(\varphi_j)_{j\in\nat_0}$ be a partition of unity as above, $\bm{w}=(w_j)_{j\in\nat_0}\in\mathcal{W}^{\alpha}_{\alpha_1,\alpha_2}(\rn)$ and $p,q \in \Plog$.
\begin{list}{}{\labelwidth1.3em\leftmargin2em}
  \item[{\upshape (i)\hfill}] The space $\Bwpqpunkt$ is defined as the collection of all $f\in \SSn$ such that
    \begin{align}
      \|f\mid \Bwpqpunkt \| &:= \Vert (w_j\,(\varphi_j \widehat{f})^{\vee} )_{j\in\nat_0}\mid \ellqp \Vert\notag
    \end{align}
    is finite.
  \item[{\upshape (ii)\hfill}] If $p^{+}, q^+< \infty$, then the space $\Bwpqpunkt$ is defined as the collection of all $f\in \SSn$ such that
    \begin{align}
      \|f \mid \Fwpqpunkt \| &:= \Vert (w_j\,(\varphi_j \widehat{f})^{\vee} )_{j\in\nat_0}\mid \ellpq \Vert\notag
    \end{align}
    is finite.
\end{list}
\end{definition}
As before, the independence of the decomposition of unity for the 2-microlocal spaces follows from the local means characterization (see \cite{Kem09} for Triebel-Lizorkin and \cite{KemV12} for Besov spaces). 

\begin{remark} These 2-microlocal weight sequences are directly connected to variable smoothness functions $s:\rn\to\mathbb{R}$ if we set
\begin{align}\label{tmlsx}
w_j(x)=2^{js(x)}.
\end{align}
If $s\in C^{\log}_{loc}(\rn)$, then $\bm{w}=(w_j(x))_{j\in\nat_0}=(2^{js(x)})_{j\in\nat_0}$ belongs to $\mathcal{W}^{\alpha}_{\alpha_1,\alpha_2}(\rn)$ with $\alpha_1=s^-$ and $\alpha_2=s^+$ and $\alpha=c_{\log}(s)$, where $c_{\log}(s)$ is the constant for $\s$ from \eqref{loc-log-holder}. That means that spaces of variable smoothness from Definition \ref{def:BFpunkt} are a special case of 2-microlocal function spaces from Definition \ref{def-Wspaces}. Both types of function spaces are very closely connected and the properties used in the proofs are either
\begin{align}\label{sxClogeigenschaft}
	2^{j|s(x)-s(y)|}\leq c\quad\text{or}\quad\frac{w_j(x)}{w_j(y)}\leq c
\end{align}
for $|x-y|\leq c\,2^{-j}$ and $j\in\nat_0$. This property follows directly either from the definition of $s\in C^{\log}_{loc}(\rn)$ or from Definition \ref{def-ad-weight}.\\
\end{remark}

\begin{theorem}
 Let $\bm{w}^0, \bm{w}^1 \in \mathcal{W}^{\alpha}_{\alpha_1,\alpha_2}(\rn)$ and $p_0, p_1, q\in \Plog$ with $q^+<\infty$. Let $p_0(x)<p_1(x)$ with $\inf_{x\in\rn}(p_1(x)-p_0(x))>0$ and
 \beq
    1< \frac{w^0_j(x)}{w^1_j(x)}=2^{j\left(\frac{n}{p_0(x)}-\frac{n}{p_1(x)}\right)}\qquad \text{for all }\, x\in\rn \text{ and } j\in\nat_0. \nonumber
 \eeq
 \begin{list}{}{\labelwidth1.3em\leftmargin2em}
  \item[{\upshape (i)\hfill}] {If $p_0^+<\infty$}, then 
    $$
    F^{w^0(\cdot)}_{p_0(\cdot),\q}\hookrightarrow B^{w^1(\cdot)}_{p_1(\cdot),p_0(\cdot)}.
    $$
  \item[{\upshape (ii)\hfill}] If $p_1^+<\infty$, then
    $$
    B^{w^0(\cdot)}_{p_0(\cdot),p_1(\cdot)}\hookrightarrow F^{w^1(\cdot)}_{p_1(\cdot),\q}.
    $$
 \end{list}
\end{theorem}
Regarding the proof, one just needs to use the corresponding Sobolev embeddings for 2-microlocal spaces (see \cite{GMN14} and \cite{AC15}) and follow exactly the same steps as before using always property \eqref{sxClogeigenschaft} for the weight sequences. 

\section{Open Problems}
We close by listing several open problems, which are connected to the
study of function spaces with variable exponents
and their embeddings.
\begin{enumerate}
\item
Give an example that $\inf_{x\in\rn}(s_0(x)-s_1(x))>0$ is really needed
and can not be replaced by $s_0(x)>s_1(x)$ for all $x\in\rn$.
This seems to be feasible when working with function space on the whole
$\rn$, but it might get more tricky, when
considering only functions with support in, say, the unit cube $[0,1]^n.$
\item It is well known, that the Triebel-Lizorkin spaces $F^s_{p,q}(\rn)$
with constant indices might depend on the chosen decomposition of unity
if $p=\infty.$ Therefore, the restriction $p^+<\infty$ seems to be quite
natural in Definition \ref{def:BFpunkt} (i).
On the other hand, there is no such trouble for spaces $F^s_{p,q}(\rn)$
with $p<\infty$ and $q=\infty.$
It would be therefore highly interesting if the Triebel Lizorkin spaces
$\Fspqpunkt$ can also be defined with $q^+=\infty$ but with still variable
$q$.
To that end one needs to show that the spaces $\Fspqpunkt$ are independent
on the resolution of unity.
\item Function spaces of Morrey type attracted recently a lot of attention
in connection with the analysis of Navier-Stokes equations \cite{GT}
and function spaces of Morrey type with variable exponents were introduced
and studied already \cite{FX}.
Is there a version of the Franke and Jawerth embedding for this scale of function spaces?
\item The Franke and Jawerth embeddings were used in \cite{DDH} to describe
the fine properties of
Besov and Triebel-Lizorkin spaces in term of the so-called envelopes.
They determine the kind and size of singularities, which the functions from
these spaces might posses.
On the other hand, function spaces of variable exponents capture very well
the local properties of functions and distributions.
The interplay of the theory of envelopes and the function spaces with
variable exponents would be therefore interesting.
\end{enumerate}

\end{document}